# Subsampling needlet coefficients on the sphere

P. BALDI[1], G. KERKYACHARIAN[2], D. MARINUCCI[1] and D. PICARD[2]

[1]*Dipartimento di Matematica, Università di Roma Tor Vergata, Italy.*
*E-mail: marinucc@mat.uniroma2.it*
[2]*Laboratoire de Probabilités et Modèles Aléatoires, Paris, France.*

In a recent paper, we analyzed the properties of a new kind of spherical wavelets (called needlets) for statistical inference procedures on spherical random fields; the investigation was mainly motivated by applications to cosmological data. In the present work, we exploit the asymptotic uncorrelation of random needlet coefficients at fixed angular distances to construct subsampling statistics evaluated on Voronoi cells on the sphere. We illustrate how such statistics can be used for isotropy tests and for bootstrap estimation of nuisance parameters, even when a single realization of the spherical random field is observed. The asymptotic theory is developed in detail in the high resolution sense.

*Keywords:* random fields; spherical wavelets; subsampling; Voronoi cells

## 1. Introduction

Strong empirical motivations have prompted rising activity in the statistical analysis of spherical random fields over the last few years, mostly in connection with the analysis of the cosmic microwave background (CMB) radiation. The first CMB maps were provided by the NASA mission *COBE* in 1993 and led to the Nobel prize for physics for J. Mather and G. Smoot in 2006. Much more refined observations were provided by the satellite mission *WMAP* in 2003/2006. Further improvements are expected in the next few years, in view of the forthcoming launch of the European Space Agency mission *Planck*. Cosmic microwave background data have been collected by many other remarkable balloon-borne experiments; a list with links to some data sets can be found on http://www.fisica.uniroma2.it/˜cosmo/.

The analysis of CMB data provides a sort of goldmine of new challenges for statistical methodology; see, for instance, Dodelson [8] for a review. Many of these challenges can be addressed by some forms of spherical wavelets; to list just a few, we recall testing for non-Gaussianity (Cabella *et al.* [6], Vielva *et al.* [29], McEwen *et al.* [17], Jin *et al.* [13]), component separation (Moudden *et al.* [18]), foreground subtraction (Hansen *et al.* [11]), point sources detection (Sanz *et al.* [25]), cross-correlation with large scale structure data (Pietrobon *et al.* [22], McEwen *et al.* [16]), searching for asymmetries/directional features in CMB (Wiaux *et al*. [32], Vielva *et al*. [30]) and many others. The importance of wavelets in this environment is easily understood: on one hand, nearly all predictions from theoretical physics for the behaviour of the CMB field are presented in Fourier space, where the orthogonality properties make many difficult problems more tractable. On the other hand, spherical CMB maps are usually observed with gaps due to the presence of foreground radiation such as the emissions by the Milky Way and other galaxies. The double localization properties of wavelets thus make them a very valuable asset for CMB data analysis. Because of this, indeed very many physical papers have attempted to use wavelets in cosmology for a variety of different problems. Many of these papers provide insightful suggestions and important experimental results; however, the focus has typically been on physical data and few have investigated the statistical foundations of the proposed procedures. In most cases, rather than proposing anew procedures suitable for spherical random fields, it has been very common to rely on tangent plane approximations to adapt to the sphere standard wavelet constructions on the plane (an exception is provided by Sanz *et al*. [25]; see also Antoine and Vandergheynst [1], Antoine *et al*. [2] and Wiaux *et al*. [31] for a nice group-theoretical construction).

Needlets are a new kind of second-generation spherical wavelets, which were introduced into functional analysis by Narcowich, Petrushev and Ward [19, 20]; they can be shown to make up a tight frame with excellent localization properties in both the real and the harmonic domains. In Baldi *et al*. [4], their properties for the analysis of random fields were established; in particular, it was shown that a major feature of random needlet coefficients is their asymptotic uncorrelation at large frequencies $j$ for any fixed angular distance (see Baldi *et al*. [3] for an analogous result on the circle). Of course, in the Gaussian and isotropic cases, this property implies that the random spherical needlets behave asymptotically as an i.i.d. array. This surprising result clearly opens the way to the implementation of several statistical procedures with an asymptotic justification. The meaning of asymptotics in this framework should be understood with great care. It should be stressed that we are considering a single observation (our universe) of a mean square continuous and isotropic random field on a fixed domain. Our asymptotic theory is then entertained in the high-resolution sense (compare Marinucci [15]), that is, as higher and higher frequency data become available for statistical analysis.

In this paper, we build on this essential feature of the random needlet coefficients to propose new statistical procedures and to provide bootstrap estimates for the asymptotic variance of existing techniques. The main idea of our work can be described as follows: We partition the sphere $\mathbb{S}^2$ into disjoint subsets with roughly the same dimension, in a sense to be made rigorous later. It is then possible to evaluate nonlinear statistics on each of these subregions separately. These subregions will constitute a tessellation of the spherical surface composed by Voronoi cells associated with a suitable $\varepsilon$-net. The geometry of these Voronoi cells on the sphere is interesting by itself and plays an important role in our results.

We show below that, in the high-resolution sense, these statistics converge asymptotically to a sequence of independent realizations with the same law as the corresponding ones evaluated on the full sphere. It is then immediate to exploit this result to obtain computationally feasible approximations of sample variance for nonlinear functionals (in Baldi *et al.* [4] estimation of the normalizing factors was left open). It is also quite straightforward to exploit our construction to derive tests for statistical isotropy on the sphere. The latter issue is indeed of great importance for cosmological data analysis: indeed the single most surprising result from the first releases of WMAP data in 2003/2006 was the apparent presence of an asymmetry in CMB radiation (see, for instance, Hansen *et al.* [11]). These findings have sparkled an impressive amount of empirical research over the last few years, but the results are still inconclusive, partly due to the lack of widely accepted statistical procedures to tackle this issue. The actual discovery of an asymmetry in cosmological radiation might call for revolutionary advances in theoretical physics, possible explanations ranging from the higher-dimension non-trivial topological structure of the observed universe to rotating solutions of Einstein field equations.

To some extent, our approach can be viewed as a simple form of subsampling in the sense of Politis *et al.* [23]. To the best of our knowledge, subsampling techniques have so far been considered mainly in the presence of the usual large sample asymptotics, the main instrument to establish asymptotic properties being some mixing properties as the observations increase. In the present circumstances, no mixing properties can be advocated, as we assume we are observing a single realization of a random field on a bounded and compact domain. Once again, then, we consider our results to be a quite surprising consequence of the peculiar properties of the needlets construction.

The structure of the paper is as follows: in Section 2 we review some basic features of the needlets and their properties in the analysis of isotropic random fields. Section 3 discusses the geometry of the sphere and, more precisely, the separation properties of Voronoi cells and cubature points inside them. In Section 4 we provide our main theoretical results, which are made possible by careful correlation inequalities that are justified in the Appendix. Section 5 discusses the results, their statistical applications and possible routes for further research.

## 2. A review on needlets

In this section, we review very briefly a few basic features of needlet construction. For more details, we refer to Narcowich, Petrushev and Ward [19] and Baldi *et al.* [4]. Let us denote by $\{Y_{lm}\}_{m=-l,\ldots,l}$ the spherical harmonics (see [28]), that form an orthonormal basis for $\mathbb{L}_2(\mathbb{S}^2, \mu)$. The following decomposition holds:

$$\mathbb{L}_2(\mathbb{S}^2, \mu) = \bigoplus_{l=0}^{\infty} H_l,$$

where the $H_l$'s are the finite dimensional spaces of $\mathbb{L}_2$ spanned by the $l$th spherical harmonics and $\mu$ is Lebesgue measure on the sphere. We define also the space of the restrictions to $\mathbb{S}^2$ of the polynomials of a degree not greater than $l$ as

$$\mathcal{K}_l = \bigoplus_{m=0}^{l} H_m.$$

Now let $\phi$ be a $C^\infty$ non-increasing function supported in $|\xi| \le 1$, such that $1 \ge \phi(\xi) \ge 0$ and $\phi(\xi) = 1$ if $|\xi| \le B^{-1}$, and define:

$$b^2(\xi) = \phi\left(\frac{\xi}{2}\right) - \phi(\xi) \ge 0$$

so that

$$\forall |\xi| \ge 1, \qquad \sum_j b^2\left(\frac{\xi}{B^j}\right) = 1, \qquad B > 1. \tag{1}$$

Define the projection operator

$$L_l(\langle x, y \rangle) = \sum_{m=-l}^{l} Y_{lm}(x)\overline{Y_{lm}(y)},$$

where we take $x, y$ on the unit sphere $\mathbb{S}^2$, $\langle x, y \rangle$ is the standard scalar product on $\mathbb{R}^3$. Let us now define the kernel

$$\Lambda_j(d(x,y)) = \sum_{[B^{j-1}] < l < [B^{j+1}]} b^2\left(\frac{l}{B^j}\right) L_l(\langle x, y \rangle),$$

$[\cdot]$ denoting integer part. Finally, we introduce *cubature points*, that is, for each $l$ we consider $\{\eta : \eta \in \mathcal{X}_l\}$, a finite subset of $\mathbb{S}^2$, and positive real numbers $\lambda_\eta > 0$ (the *cubature weights*) indexed by the elements $\eta$ of $\mathcal{X}_l$, such that

$$\forall f \in \mathcal{K}_l, \qquad \int f \, d\mu = \sum_{\eta \in \mathcal{X}_l} \lambda_\eta f(\eta). \tag{2}$$

It can be proved that, if $\mathcal{X}_l$ is a maximal $\varepsilon$-net (see the definition below) with $\varepsilon \sim B^{-j}$ for some $j$, then it constitutes a set of cubature points (see [19]). This will be our choice in the sequel. We will write $\mathcal{Z}_j = \mathcal{X}_{2[B^{j+1}]}$, $\mathcal{Z}_j = \{\xi_{j1}, \ldots, \xi_{jk}, \ldots\}$ and $\lambda_{jk} = \lambda_\eta$ in order to stick with the usual wavelet notation.

We are now ready to introduce needlets, which are given by

$$\psi_{jk}(x) := \sum_{[B^{j-1}] < l < [B^{j+1}]} \sqrt{\lambda_{jk}}\, b\left(\frac{l}{B^j}\right) L_l(\langle x, \xi_{jk} \rangle), \qquad j = 1, 2, \ldots, \xi_{jk} \in \mathcal{Z}_j.$$

Needlets enjoy a number of very important properties that are proved elsewhere and that we summarize in the following two propositions. Write as usual

$$\langle f, \psi_{jk} \rangle_{L^2(\mathbb{S}^2)} := \int_{\mathbb{S}^2} f(x) \psi_{jk}(x) \, d\mu;$$

we have the following proposition.

**Proposition 1 ([19]).** (a) (*Reconstruction*) *The family* $(\psi_{j\eta})_{j\in\mathbb{N},\eta\in\mathcal{Z}_j}$ *is a tight frame, hence*

$$f(x) = \sum_j \sum_{\xi_{jk}\in\mathcal{Z}_j} \langle f, \psi_{jk}\rangle_{L^2(\mathbb{S}^2)} \psi_{jk}(x). \tag{3}$$

(b) *(Localization) For any positive integer $M$ there exists a constant $c_M$ such that:*

$$|\psi_{jk}(x)| \le \frac{c_M B^j}{[1+B^j \mathrm{d}(x,\xi_{jk})]^M}, \tag{4}$$

*where $d(x,y) = \arccos(\langle x,y\rangle)$ denotes the usual geodesic distance on the sphere.*

Tight frames in some sense are very close to orthonormal bases, as shown by (3); see [12] for further references and discussion. The good localization properties highlighted in (4) are very useful when dealing with missing observations or with statistics evaluated on subsets of the sphere (see below).

Let us now focus on zero-mean, mean square continuous and isotropic random fields. As is well known, these can be expressed by means of the following spectral decomposition, which holds in the mean square sense:

$$T(x) = \sum_{l=1}^{\infty} T_l(x) = \sum_{lm} a_{lm} Y_{lm}(x),$$

where $\{a_{lm}\}_{l,m}$, $m=1,\ldots,l$ is a triangular array of zero-mean, orthogonal, complex-valued random variables with variance $\mathbb{E}|a_{lm}|^2 = C_l$, the angular power spectrum of the random field. For $m<0$ we have $a_{lm} = (-1)^m \overline{a_{l-m}}$, whereas $a_{l0}$ is real. Our random needlet coefficients are then easily seen to be given by

$$\beta_{jk} := \int_{\mathbb{S}^2} T(x)\psi_{jk}(x)\,\mathrm{d}x = \sqrt{\lambda_{jk}} \sum_{[B^{j-1}]<l<[B^{j+1}]} b\left(\frac{l}{B^j}\right) T_l(\xi_{jk}), \qquad \xi_{jk}\in\mathcal{Z}_j.$$

Remark that in the previous sum only a finite number of non-vanishing terms appear. Hence we obtain

$$\mathbb{E}\beta_{jk}\beta_{jk'} = \sqrt{\lambda_{jk}\lambda_{jk'}} \sum_{l\ge 0} b^2\left(\frac{l}{B^j}\right) C_l L_l(\langle \xi_{j,k},\xi_{j,k'}\rangle)$$

and

$$\mathrm{corr}(\beta_{j,k},\beta_{j,k'}) = \frac{\sum_l b^2(l/B^j) C_l L_l(\langle \xi_{j,k},\xi_{j,k'}\rangle)}{\sum_l b^2(l/B^j) C_l L_l(1)} = \frac{\sum_l b^2(l/B^j) C_l L_l(\langle \xi_{j,k},\xi_{j,k'}\rangle)}{\sum_l b^2(l/B^j) C_l (2l+1)/(4\pi)}.$$

In particular $\mathbb{E}[\beta_{jk}^2]$ is equal to

$$\lambda_{jk} \underbrace{\sum_{l\geq 0} b^2\left(\frac{l}{B^j}\right) C_l L_l(1)}_{:=\gamma_j}. \tag{5}$$

The following assumption is a mild condition on the regularity of the angular power spectrum.

***Condition 1.*** *Let us assume that $\exists \alpha > 2$ and an integer $M \geq 3$ such that*

$$C_l = l^{-\alpha} g_j\left(\frac{l}{B^j}\right), \qquad \text{for } B^{j-1} < l < B^{j+1}$$

*and*

$$\forall r = 0, 1, \ldots, M, \qquad \sup_j \sup_{1/B < u < B} |g_j^{(r)}(u)| \leq C_r < \infty.$$

***Remark 2.*** A typical example that fulfills Condition 1 is

$$C_l = G(l) l^{-\alpha}, \qquad \sup_{u \geq 1/B} u^r |G^{(r)}(u)| \leq C_r < \infty, \qquad r = 0, 1, \ldots, M.$$

The following proposition provides the basic correlation inequality for the analysis of the statistical procedures in the sequel.

**Proposition 3 ([4]).** *Under Condition 1 we have, for each positive integer $M$,*

$$|\operatorname{corr}(\beta_{j,k}, \beta_{j,k'})| \leq \frac{C_M B^{2j}}{(1 + B^j d(\xi_{j,k}, \xi_{j,k'}))^M}. \tag{6}$$

***Remark 4.*** The previous inequality is the basic ingredient for many of the results to follow, together with the geometric analysis in the next section. Inequality (6) entails that needlet coefficients are asymptotically (in the high-resolution sense) uncorrelated at any fixed angular distance. This makes consistent inference procedures viable, even in the presence of observations on a single realization of a random field.

## 3. The geometry of Voronoi cells on the sphere

In this section, we establish some properties of Voronoi cells on the sphere. These properties will be instrumental for the correlation inequalities we shall need in our main arguments. Let us first introduce some notation.

We start by defining the standard (open and closed) balls in $\mathbb{S}^2$ as

$$B(a, \alpha) = \{x, \ d(a, x) \leq \alpha\}, \qquad B^\circ(a, \alpha) = \{x, \ d(a, x) < \alpha\}.$$

Also, if $A \subset \mathbb{S}^2$ we denote by $|A|$ the two-dimensional measure of $A$. Of course, for the full sphere $|\mathbb{S}^2| = 4\pi$. Now let $\varepsilon > 0$ and $x_1, \ldots, x_N \in \mathbb{S}^2$ such that

$$\forall i \neq j, \qquad d(x_i, x_j) > \varepsilon$$

and the set $\{x_1, \ldots, x_N\} = \Xi_\varepsilon$ is maximal for this property, that is,

$$\forall x \in \mathbb{S}^2, \qquad d(x, \Xi_\varepsilon) \leq \varepsilon \quad \text{and} \quad \forall i \neq j, \qquad B\left(x_i, \frac{\varepsilon}{2}\right) \cap B\left(x_j, \frac{\varepsilon}{2}\right) = \varnothing.$$

We call such a set $\Xi_\varepsilon$ a *maximal $\varepsilon$-net*.

Of course, if $\varepsilon = \pi$, obviously $N = 1$, and as soon as $\varepsilon < \pi$, $N \geq 2$. For instance, if $\frac{\pi}{2} \leq \varepsilon < \pi$, we can take $N = 2$ and $\xi_1$ and $\xi_2$ are two opposite poles; however, $N = 3$ is possible taking three points on a geodesic circle at distance $\frac{2\pi}{3}$. The next lemma is a simple but useful result.

**Lemma 5.** *Let $\{x_1, \ldots, x_N\} = \Xi_\varepsilon$ be a maximal $\varepsilon$-net. Then*

$$\frac{4}{\varepsilon^2} \leq N \leq \frac{4}{\varepsilon^2}\pi^2. \tag{7}$$

**Proof.** Let us first recall that $\forall a \in \mathbb{S}^2$, $\forall 0 < \eta \leq \pi$,

$$|B(a, \eta)| = 2\pi \int_0^\eta \sin\theta \, d\theta = 2\pi(1 - \cos\eta) = 4\pi \sin^2(\eta/2) \qquad (\sim \pi\eta^2 \text{ for } \eta \to 0).$$

More precisely, if $0 < \eta \leq 2\alpha \leq \pi$, as $\frac{1}{\alpha}\sin\alpha \geq \frac{2}{\pi}$,

$$\eta^2 \frac{4}{\pi} \leq \pi\eta^2 \left(\frac{\sin\alpha}{\alpha}\right)^2 \leq |B(a,\eta)| \leq \pi\eta^2. \tag{8}$$

As

$$\bigcup_{x_i \in \Xi_\varepsilon} B(x_i, \varepsilon) = \mathbb{S}^2,$$

we have

$$4\pi = |\mathbb{S}^2| \leq \sum_{x_i \in \Xi_\varepsilon} |B(x_i, \varepsilon)| \leq N\pi\varepsilon^2$$

and, as the $B(x_i, \frac{\varepsilon}{2})$ are disjoint

$$N\left(\frac{\varepsilon}{2}\right)^2 \frac{4}{\pi} \leq \sum_{x_i \in \Xi_\varepsilon} \left|B\left(x_i, \frac{\varepsilon}{2}\right)\right| \leq 4\pi. \qquad \square$$

**Remark 6.** Actually, it is easy to see that

$$\frac{1}{\sin^2 \varepsilon/2} \leq N \leq \frac{1}{\sin^2 \varepsilon/4},$$

so that

$$\frac{1}{4} \leq \frac{1}{4\cos^2 \varepsilon/4} \leq N \sin^2 \frac{\varepsilon}{4} \leq 1.$$

We can now introduce, for all $x_i \in \Xi_\varepsilon$, the Voronoi cells:

**Definition 7.** *Let $\Xi_\varepsilon$ be a maximal $\varepsilon$-net. We define the associate family of Voronoi cells:*

$$S(x_i) = \{x \in \mathbb{S}^2, \forall j \neq i, d(x, x_i) \leq d(x, x_j)\}.$$

Clearly

$$B\left(x_i, \frac{\varepsilon}{2}\right) \subset S(x_i) \subset B(x_i, \varepsilon).$$

Now observe that $\forall a \neq b$ and the set $\{x \in \mathbb{S}, d(x, a) = d(x, b)\}$ is a geodesic circle on the sphere. Let us denote this circle by $C_{a,b}$. The set: $\{x \in \mathbb{S}^2, d(x, a) \leq d(x, b)\} = D(a, b)$ is the half-sphere defined by this geodesic circle containing $a$. We have

$$S(x_i) = \bigcap_{j \neq i} D(x_i, x_j),$$

so $S(x_i)$ is the intersection of a finite number of half-spheres. Actually,

$$S(x_i) = \bigcap_{j \neq i, d(x_j, x_i) \leq 2\varepsilon} D(x_i, x_j)$$

and, more precisely,

$$S(x_i) \cap S(x_j) \neq \varnothing \quad \Longrightarrow \quad d(x_i, x_j) \leq 2\varepsilon.$$

The following lemma provides a bound on the number of points of a maximal $\varepsilon$-net that lie within a distance $2\varepsilon$ from a fixed center.

**Lemma 8.**

$$\operatorname{Card}\{j \neq i, d(x_j, x_i) \leq 2\varepsilon\} \leq 6\pi^2.$$

**Proof.** Let $K(\varepsilon, x_i) = \{x_j \neq x_i, d(x_j, x_i) \leq 2\varepsilon\}$. Clearly, as

$$K(\varepsilon, x_i) \subset B(x_i, 2\varepsilon) \setminus B(x_i, \varepsilon),$$

we have, if $x_j \in K(\varepsilon, x_i)$, $B(x_j, \frac{\varepsilon}{2}) \subset B(x_i, \frac{5\varepsilon}{2}) \setminus B(x_i, \frac{\varepsilon}{2})$. Hence, since the balls $B(x_j, \frac{\varepsilon}{2})$ are disjoint,

$$\sum_{\{j \neq i, d(x_j, x_i) \leq 2\varepsilon\}} \left| B\left(x_j, \frac{\varepsilon}{2}\right) \right| \leq \left| B\left(x_i, \frac{5\varepsilon}{2}\right) \setminus B\left(x_i, \frac{\varepsilon}{2}\right) \right|. \tag{9}$$

For $0 < \mu < \eta \leq \pi$, we have

$$|B(a, \eta) \setminus B(a, \mu)| = 2\pi \int_\mu^\eta \sin\theta \, d\theta = 2\pi(\cos\mu - \cos\eta)$$

$$= 4\pi \sin\left(\frac{\eta - \mu}{2}\right) \sin\left(\frac{\eta + \mu}{2}\right) \leq \pi(\eta^2 - \mu^2).$$

Hence

$$\eta \leq \frac{\pi}{2} \implies \frac{4}{\pi}(\eta^2 - \mu^2) \leq |B(a, \eta) \setminus B(a, \mu)| \leq \pi(\eta^2 - \mu^2). \tag{10}$$

Now, if $\eta \to 0$ and $0 \leq \mu < \eta$,

$$|B(a, \eta) \setminus B(a, \mu)| \sim \pi(\eta^2 - \mu^2). \tag{11}$$

So, by (9) and (8),

$$\text{card}(K(\varepsilon, x_i))\left(\frac{\varepsilon}{2}\right)^2 \frac{4}{\pi} \leq \pi\left(\left(\frac{5\varepsilon}{2}\right)^2 - \left(\frac{\varepsilon}{2}\right)^2\right),$$

which implies

$$\text{card}(K(\varepsilon, x_i)) \leq 6\pi^2,$$

as required. □

Of course

$$\bigcup_{x_i \in \Xi_\varepsilon} S(x_i) = \bigcup_{x_i \in \Xi_\varepsilon} B(x_i, \varepsilon) = \mathbb{S}^2.$$

It is obvious that if $x, y \in S(x_i)$ then the portion of geodesic circle joining $x$ to $y$ is inside $S(x_i)$.

**Definition 9.** *Two Voronoi cells $S(x_i)$ and $S(x_j)$ with $i \neq j$ are said to be adjacent if*

$$S(x_i) \cap S(x_j) \neq \varnothing.$$

**Remark 10.** *From Lemma 8 it is clear that there are at most $6\pi^2$ adjacent cells to any given cell.*

With the previous results at hand, we have now the ingredients to define our statistics on suitable subsets of $\mathbb{S}^2$, as explained in the following section.

## 4. Asymptotics for needlet functionals

In Baldi *et al.* [4] we discussed the asymptotics of several statistics of the form

$$\Gamma_j^{(q)} := \frac{1}{A_j} \sum_k H_q(\widehat{\beta}_{jk}), \qquad \text{where } \widehat{\beta}_{jk} := \frac{\beta_{jk}}{\sqrt{\mathbb{E}\beta_{jk}^2}}.$$

Here $A_j$ denotes the number of cubature points and $H_q$ as usual are the Hermite polynomials of order $q$, (see, e.g., Surgailis [27]):

$$H_q(u) = (-1)^q e^{u^2/2} \frac{d^q}{du^q} e^{-u^2/2},$$

that is, $H_1(u) = u$, $H_2(u) = u^2 - 1$, $H_3(u) = u^3 - 3u, \ldots$. It was shown that $\mathbb{E}\Gamma_j^{(q)} = 0$ and, as $j \to \infty$

$$\frac{\Gamma_j^{(q)}}{\sqrt{\Sigma_j^{(q)}}} \to_d N(0,1), \qquad \text{where } \Sigma_j^{(q)} := \mathbb{E}[(\Gamma_j^{(q)})^2].$$

Indeed it is also possible to show stronger results; for instance, the joint convergence over different values of $q$.

Our purpose here is to provide estimates for the normalizing variance $\Sigma_j^{(q)}$. Our idea is to "subsample" our statistics by evaluating them on distinct Voronoi cells. By means of the asymptotic properties of needlet coefficients, we shall show that each cell will provide asymptotically an independent realization of the same limiting distribution as the statistic on the whole sphere. It will then be immediate to exploit these results to estimate by bootstrap our limiting distribution. More precisely, fix $B > 1$, let $j, r$ be two non-negative integers such that $B^{-j} < B^{-r}/4$ and write $\Xi_{\pi B^{-j}} := \{x_{j,1}, \ldots, x_{j,A_j}\}, \Xi_{\pi B^{-r}} := \{x_{r,1}, \ldots, x_{r,A_r}\}$ for corresponding sequences of maximal $\varepsilon$-nets with $\varepsilon = \pi B^{-j}$ and $= \pi B^{-r}$, respectively. Let us also define

$$N_{a;rj} := \operatorname{Card}(S(x_{r,a}) \cap \Xi_{\pi B^{-j}}), \qquad A_r = \operatorname{Card}\{\Xi_{\pi B^{-r}}\},$$

where $S(x_{r,a})$ are the Voronoi cells of the $\pi B^{-r}$-net $\Xi_{\pi B^{-r}}$. In other words, $N_{ra;j}$ represents the number of points in the $B^{-j}$-net that fall inside the Voronoi cell of $\Xi_{\pi B^{-r}}$ around the point $x_{r,a}$ (note that $x_{r,a}$ need not belong to $\Xi_{\pi B^{-j}}$). $A_r$ denotes the cardinality of such Voronoi cells. In particular, $\sum_{a=1}^{A_r} N_{a;rj} = A_j$.

We focus on the triangular array

$$\Gamma_{a;rj} := \frac{1}{\sqrt{N_{a;rj}}} \sum_{k \in \Xi_{\pi B^{-j}} \cap S(x_{r,a})} H_q(\widehat{\beta}_{jk}), \qquad a = 1, 2, \ldots, A_r, \tag{12}$$

for a fixed positive integer $q$. We define $\sigma_j^2 := \text{Var}\{H_q(\beta_{jk})\}$, the variance of each summand in (12), and $\Sigma_{a;rj} := \text{Var}\{\Gamma_{a;rj}\}$ the variance of the whole sum. In the sequel, we drop the index $q$ whenever this causes no ambiguity. In particular, for $r = 0$ the unique Voronoi cell is $\mathbb{S}^2$ itself and we write $\Sigma_j := \Sigma_{a;0j}$. We prove the following:

**Theorem 11.** *Assume that the field is Gaussian and that Condition 1 holds.*

(a) *As $j \to \infty$*

$$\lim_{j \to \infty} \frac{\widehat{\sigma}_j^2}{\sigma_j^2} = 1, \qquad \text{in probability,}$$

*where*

$$\widehat{\sigma}_j^2 := \frac{1}{A_j} \sum_{k \in \Xi_{\pi B^{-j}}} H_q(\widehat{\beta}_{jk})^2.$$

(b) *For $r, j$ such that*

$$\frac{1}{r} + \frac{r}{j} \to 0, \tag{13}$$

$\Sigma_j$ *is consistently estimated by*

$$\widehat{\Sigma}_j = \frac{1}{A_r} \sum_{a=1}^{A_r} \Gamma_{a;rj}^2, \tag{14}$$

*that is,*

$$\lim_{j \to \infty} \frac{\widehat{\Sigma}_j}{\Sigma_j} = 1, \qquad \text{in probability.}$$

(c) *For $r$ (and hence $A_r$) fixed and $j \to \infty$, we have*

$$\left\{\frac{\Gamma_{1;rj}}{\sqrt{\Sigma_j}}, \ldots, \frac{\Gamma_{A_r;rj}}{\sqrt{\Sigma_j}}\right\}' \to_d N_{A_r}(0, I), \qquad \text{as } j \to \infty. \tag{15}$$

(d) *As $r \to \infty$, and $\frac{j}{r} \to \infty$, we have*

$$\left\{\frac{\Gamma_{ra;j}}{\sqrt{\Sigma_{rj}}}\right\}_{a=1,\ldots,A_r} \to \{Z_a\}_{a=1,\ldots,A_r}, \qquad Z_a \stackrel{d}{=} NID(0,1), \tag{16}$$

$$\left\{\frac{\Gamma_{ra;j}}{\sqrt{\widehat{\Sigma}_{rj}}}\right\}_{a=1,\ldots,A_r} \to \{Z_a\}_{a=1,\ldots,A_r}, \qquad Z_a \stackrel{d}{=} NID(0,1), \tag{17}$$

*the convergence taking place in the space of sequences $\ell^\infty$ endowed with the standard topology.*

***Remark 12.*** We note here that (a) and (b) of Theorem 11 provide some sort of bootstrap/subsampling approximation for the sample variance for the statistics

$$\Gamma_j = \frac{1}{A_j} \sum_{k \in \Xi_{\pi B^{-j}}} H_q(\widehat{\beta}_{jk}),$$

thus eliminating the need for parametric assumptions on nuisance parameters or lengthy Monte Carlo simulations. In particular, for $q = 1$, point (a) yields a consistent estimate of the normalizing factor $\gamma_j$ defined in (5).

Condition (13) of Theorem 11 is clearly reminiscent of standard bandwidth assumptions in nonparametric estimation. It ensures that the number of pixels on which needlet coefficients can be evaluated grows more rapidly than the number of Voronoi cells or, in other words, the number of observations within each Voronoi cell goes to infinity. Points (c) and (d) entail that statistics evaluated on (possibly adjacent) Voronoi cells are asymptotically i.i.d.

In the sequel, we take $r$ to be an increasing function of $j$, but we refrain from writing $r_j$ to simplify notation.

The proof of Theorem 11 requires a very careful bound on cross-correlation of our statistics over different Voronoi cells. Such bounds rely very much on the following result, the proof of which is quite delicate and collected in the Appendix.

**Proposition 13.** *Let $0 < \varepsilon \leq \delta/4$ and let $\Xi_\delta$ be a maximal $\delta$-net and $\Xi_\varepsilon$ be a maximal $\varepsilon$-net. For $x_a \in \Xi_\delta$, let $S(x_a)$ denote the Voronoi cells centered at $x_a$ for the $\delta$-net $\Xi_\delta$. Let $N_{a;\delta\varepsilon} = \mathrm{Card}(S(x_a) \cap \Xi_\varepsilon)$ and assume that $\varepsilon = \pi B^{-j}$. Then, if $x_a, x_b \in \Xi_\delta$, $x_a \neq x_b$ and $M \geq 3$, there exists $C > 0$ verifying the following property:*

$$\frac{1}{\sqrt{N_{a;\delta\varepsilon} N_{b;\delta\varepsilon}}} \sum_{v \in \Xi_\varepsilon \cap S(x_a)} \sum_{u \in \Xi_\varepsilon \cap S(x_b)} \frac{1}{(1 + B^j d(u,v))^M} \leq C \frac{\varepsilon}{\delta} \left( 1 + 1_{\{M=3\}} \log \frac{\delta}{\varepsilon} \right).$$

**Proof of Theorem 11.** From Proposition 13 it is immediate to get the following inequality, which we shall exploit several times in our arguments below. For some $C > 0$, we have

$$\begin{aligned}\frac{1}{\sqrt{N_{a;\delta\varepsilon} N_{b;\delta\varepsilon}}} \sum_{v \in \Xi_\varepsilon \cap S(x_a)} \sum_{u \in \Xi_\varepsilon \cap S(x_b)} \frac{1}{(1 + B^j d(u,v))^M} \\ \leq C(j-r) B^{-(j-r)} \log B.\end{aligned} \quad (18)$$

(a) The proof of (a) is very similar to the proof of b), indeed slightly simpler, and hence omitted.

(b) It is readily checked that $\lim_{j \to \infty} \mathbb{E}\{\widehat{\Sigma}_j / \Sigma_j\} = 1$. The idea of our argument is then to establish

$$\limsup_{j \to \infty} \mathrm{Var}\left\{ \frac{\Gamma_{ra;j}^2}{\Sigma_j} \right\} = \mathrm{O}(1), \quad (19)$$

$$\limsup_{j \to \infty} \max_{a,b} \left| \operatorname{Cov}\left\{ \frac{\Gamma_{ra;j}^2}{\Sigma_j}, \frac{\Gamma_{rb;j}^2}{\Sigma_j} \right\} \right| = 0 \tag{20}$$

and then to use these results to conclude the proof by noting that

$$\limsup_{j \to \infty} \operatorname{Var}\left\{ \frac{\widehat{\Sigma}_j}{\Sigma_j} \right\} = \limsup_{j \to \infty} \frac{1}{A_r^2} \sum_{a,b=1}^{A_r} \operatorname{Cov}\left\{ \frac{\Gamma_{ra;j}^2}{\Sigma_j}, \frac{\Gamma_{rb;j}^2}{\Sigma_j} \right\}$$

$$= \limsup_{j \to \infty} \frac{1}{A_r^2} \sum_{a=1}^{A_r} \operatorname{Var}\left\{ \frac{\Gamma_{ra;j}^2}{\Sigma_j} \right\} + \operatorname{O}\left( \frac{1}{A_r^2} \sum_{a \neq b=1}^{A_r} \operatorname{Cov}\left\{ \frac{\Gamma_{ra;j}^2}{\Sigma_j}, \frac{\Gamma_{rb;j}^2}{\Sigma_j} \right\} \right) = 0.$$

To make this argument rigorous, we start by noting that

$$\operatorname{Cov}\left\{ \frac{\Gamma_{ra;j}^2}{\Sigma_j}, \frac{\Gamma_{rb;j}^2}{\Sigma_j} t \right\}$$

$$= \frac{1}{N_{ra;j} N_{rb;j}} \tag{21}$$

$$\times \sum_{\substack{k_1, k_2 \in \Xi_{\pi B^{-j}} \cap S(x_{ra}) \\ k_3, k_4 \in \Xi_{\pi B^{-j}} \cap S(x_{rb})}} \frac{\mathbb{E}\{H_q(\beta_{jk_1}) H_q(\beta_{jk_2}) H_q(\beta_{jk_3}) H_q(\beta_{jk_4})\}}{\Sigma_j^2}.$$

Let $\rho_j(k_1, k_2) = \operatorname{corr}(\beta_{jk_1}, \beta_{jk_2}) = \operatorname{corr}(\widehat{\beta}_{jk_1}, \widehat{\beta}_{jk_2})$. By the diagram formula (see, e.g., Surgailis [27]) and with standard manipulations it can be shown that (21) is bounded by

$$\frac{Kq!}{N_{ra;j} N_{rb;j}} \sum_{\substack{k_1, k_2 \in \Xi_{\pi B^{-j}} \cap S(x_{ra}) \\ k_3, k_4 \in \Xi_{\pi B^{-j}} \cap S(x_{rb})}} |\rho_j^{q_1}(k_1, k_2) \rho_j^{q_2}(k_1, k_3)$$

$$\times \rho_j^{q_3}(k_1, k_4) \rho_j^{q_4}(k_2, k_3) \rho_j^{q_5}(k_2, k_4) \rho_j^{q_6}(k_3, k_4)| \tag{22}$$

for some $K > 0$ and for all choices of non-negative integers $q_1, \ldots, q_6$ such that $q_1 + \cdots + q_6 = 2q$ and $(q_2 \vee q_3), (q_4 \vee q_5)$ are strictly positive. Rearranging terms and recalling that $|\rho_j(\cdot, \cdot)| \leq 1$, we obtain

$$(22) \leq \frac{K}{N_{ra;j} N_{rb;j}} \sum_{\substack{k_1, k_2 \in \Xi_{\pi B^{-j}} \cap S(x_{ra}) \\ k_3, k_4 \in \Xi_{\pi B^{-j}} \cap S(x_{rb})}} |\rho_j(k_1, k_3) \rho_j(k_2, k_4)|$$

$$= \left\{ \sqrt{\frac{K}{N_{ra;j} N_{rb;j}}} \sum_{\substack{k_1 \in \Xi_{\pi B^{-j}} \cap S(x_{ra}) \\ k_3 \in \Xi_{\pi B^{-j}} \cap S(x_{rb})}} |\rho_j(k_1, k_3)| \right\}^2$$

$$\leq \left\{\sqrt{\frac{K}{N_{ra;j}N_{rb;j}}} \sum_{\substack{k_1 \in \Xi_{\pi B^{-j}} \cap S(x_{ra}) \\ k_3 \in \Xi_{\pi B^{-j}} \cap S(x_{rb})}} \frac{C_M B^{2j}}{(1+B^j d(\xi_{j,k_1},\xi_{j,k_3}))^M}\right\}^2$$

$$\leq C(j-r)^2 B^{-2(j-r)} \log^2 B,$$

in view of Proposition 13 and (18). Thus (20) is established. The proof of (19) is similar. It follows that

$$\lim_{j\to\infty} \operatorname{Var}\left\{\frac{\widehat{\Sigma}_j}{\Sigma_j}\right\} = \lim_{j\to\infty} \operatorname{Var}\left\{\frac{1}{A_r}\sum_{a=1}^{A_r}\frac{\Gamma_{ra;j}^2}{\Sigma_j}\right\}$$

$$= \lim_{j\to\infty} \frac{1}{A_r^2}\sum_{a=1}^{A_r} \operatorname{Var}\left\{\frac{\Gamma_{ra;j}^2}{\Sigma_j}\right\} + \lim_{j\to\infty} \frac{1}{A_r^2}\sum_{a,b=1}^{A_r} \operatorname{Cov}\left\{\frac{\Gamma_{ra;j}^2}{\Sigma_j}, \frac{\Gamma_{rb;j}^2}{\Sigma_j}\right\}$$

$$\leq \operatorname{O}\left(\frac{1}{A_r}\right) + \limsup_{j\to\infty} \max_{1\leq a,b\leq A_r}\left|\operatorname{Cov}\left\{\frac{\Gamma_{ra;j}^2}{\Sigma_j}, \frac{\Gamma_{rb;j}^2}{\Sigma_j}\right\}\right| = \operatorname{o}(1), \qquad \text{as } j\to\infty,$$

which completes the proof of (b).

(c) Here we recall that $A_j = A$ is fixed, that is, we are focusing on a finite number of subsets of the sphere. In these circumstances, the argument is very similar to the proof of Theorem 9 in Baldi et al. [4]. We use the Cramér–Wold device and thus focus on the linear combination

$$\frac{1}{\sqrt{\Sigma_j}}\sum_{a=1}^{A_r} w_{ja}\Gamma_{ja} \to_d N\left(0, \sum_{a=1}^{A_r} w_{ja}^2\right),$$

where $w_{ja} \in \mathbb{R}$ for all $j=1,2,\ldots, a=1,\ldots,A_j$. It is obvious that

$$\mathbb{E}\left\{\frac{1}{\sqrt{\Sigma_j}}\sum_{a=1}^{A_r} w_{ja}\Gamma_{ja}\right\} = 0.$$

On the other hand

$$\operatorname{Var}\left\{\frac{1}{\sqrt{\sum_{a=1}^{A_r} w_{ja}^2 \Sigma_j}}\sum_{a=1}^{A_r} w_{ja}\Gamma_{ja}\right\}$$

$$= \mathbb{E}\left\{\frac{1}{\sqrt{\sum_{a=1}^{A_r} w_{ja}^2 \Sigma_j}}\sum_{a=1}^{A_r} w_{ja}\Gamma_{ja}\right\}^2$$

$$= 1 + \frac{1}{\sum_{a=1}^{A_r} w_{ja}^2 \Sigma_j}\sum_{a\neq b} w_{ja}w_{jb}\mathbb{E}\{\Gamma_{ja}\Gamma_{jb}\} \to 1,$$

as $j \to \infty$, in view of (20). Also, let us denote by $\text{cum}_p(X)$ the cumulant of order $p$ of the random variable $X$; to establish the central limit theorem, we resort to a recent result by Nualart and Peccati [21], where it is proved that convergence to zero of the fourth-order cumulant is a sufficient condition for asymptotic Gaussianity in the Gaussian subordinated case we are considering here (see also DeJong [7] for a similar approach with multilinear forms in i.i.d. sequences). Again by the diagram formula for cumulants, we obtain easily for $p = 4$,

$$\text{cum}_4\left\{\sum_{a=1}^{A_r} w_{ja} \frac{\Gamma_{ja}}{\sqrt{\Sigma_j}}\right\}$$

$$= \sum_{a,b,c,d=1}^{A_r} w_{ja} w_{jb} w_{jc} w_{jd} \,\text{cum}\left\{\frac{\Gamma_{ja}}{\sqrt{\Sigma_j}}, \frac{\Gamma_{jb}}{\sqrt{\Sigma_j}}, \frac{\Gamma_{jc}}{\sqrt{\Sigma_j}}, \frac{\Gamma_{jd}}{\sqrt{\Sigma_j}}\right\}$$

$$\leq \frac{K A_r^4}{\sqrt{N_{ra;j} N_{rb;j} N_{rc;j} N_{rd;j}}}$$

$$\times \sum_{\substack{k_1 \in \Xi_{\pi B^{-j}} \cap S(x_{ra}) \\ \cdots \\ k_4 \in \Xi_{\pi B^{-j}} \cap S(x_{rd})}} \text{cum}\left\{\frac{H_q(\beta_{jk_1})}{\sqrt{\Sigma_j}}, \frac{H_q(\beta_{jk_2})}{\sqrt{\Sigma_j}}, \frac{H_q(\beta_{jk_3})}{\sqrt{\Sigma_j}}, \frac{H_q(\beta_{jk_4})}{\sqrt{\Sigma_j}}\right\}$$

$$\leq \frac{K A_r^4}{\sqrt{N_{ra;j} N_{rb;j} N_{rc;j} N_{rd;j}}}$$

$$\times \sum_{\substack{k_1 \in \Xi_{\pi B^{-j}} \cap S(x_{ra}) \\ \cdots \\ k_4 \in \Xi_{\pi B^{-j}} \cap S(x_{rd})}} |\rho_j(k_1, k_2)||\rho_j(k_2, k_3)||\rho_j(k_3, k_4)||\rho_j(k_4, k_1)|$$

$$\leq K A_r^4 \left\{\frac{1}{\sqrt{N_{ra;j} N_{rb;j}}} \sum_{\substack{k_1 \in \Xi_{\pi B^{-j}} \cap S(x_{ra}) \\ k_2 \in \Xi_{\pi B^{-j}} \cap S(x_{rb})}} \frac{C_M B^{2j}}{(1 + B^j d(\xi_{j,k_1}, \xi_{j,k_2}))^M}\right\}$$

$$\times \left\{\frac{1}{\sqrt{N_{rc;j} N_{rd;j}}} \sum_{\substack{k_3 \in \Xi_{\pi B^{-j}} \cap S(x_{rc}) \\ k_4 \in \Xi_{\pi B^{-j}} \cap S(x_{rd})}} \frac{C_M B^{2j}}{(1 + B^j d(\xi_{j,k_3}, \xi_{j,k_4}))^M}\right\}$$

$$= o(1), \qquad \text{as } j \to \infty,$$

again in view of Proposition 13 and (18). The central limit theorem is then established.

(d) Fix a finite subset $D \subset \mathbb{N}$ of cardinality $A$ and label its elements $a = 1, \ldots, A$. It is obvious that as $j \to \infty$

$$\left\{\frac{\Gamma_{rk;j}}{\sqrt{\Sigma_j}}\right\}_{k \in D} = \left\{\frac{\Gamma_{ra;j}}{\sqrt{\Sigma_j}}\right\}_{a=1,\ldots,A} \to_d Z_a \stackrel{d}{=} N(0,1),$$

by exactly the same argument as in (c). Because $D$ is arbitrary, we have thus established convergence in the finite-dimensional distributions sense of the sequence $\{\Gamma_{ra;j}/\sqrt{\Sigma_j}\}_{a\in\mathbb{N}}$ to the sequence $\{Z_a\}_{a\in\mathbb{N}}$. It is a standard result that in $\mathbb{R}^\infty$ the finite-dimensional sets are a determining class (see Billingsley [5], page 19), so that weak convergence is actually equivalent to convergence of the finite-dimensional distributions. The argument for $\{\Gamma_{ra;j}/\sqrt{\widehat{\Sigma}_j}\}_{a\in\mathbb{N}}$ is entirely analogous by means of Slutzky's lemma.

□

## 5. Discussion and statistical applications

The results in the previous section lend themselves to several applications for the statistical analysis of spherical random fields, in particular for CMB data. Here we briefly discuss some examples; we do not provide a complete discussion or application to real data, as we prefer to defer this more detailed investigation to future works. Refer also to [9, 22] for numerical evidence and applications of needlets to CMB data.

An immediate application of the results in Theorem 11 concerns studentization. Indeed, in Baldi *et al.* [4], several statistics were proposed, based on needlet functionals, to test for goodness of fit, Gaussianity and isotropy of CMB data. We remark that these three issues are very widely studied in the huge satellite experiment collaborations on CMB, such as WMAP and Planck. These statistics could all be expressed as linear combinations of $\{\Gamma_{ra;j}\}$ in (12), and asymptotic Gaussianity was established under the null, provided the normalizing variance could be taken to be known. The results of the previous section can then be immediately exploited to provide estimates of the limiting variances, thus making feasible testing procedures with a standard asymptotic distribution. In particular, Theorem 11, part (d) shows how nonlinear statistics can be studentized by the estimates $\widehat{\Sigma}_j$ given in (14) without affecting the limiting distribution.

Other possible applications of our results relate to testing for isotropy of spherical random fields data. Indeed, in cosmological applications an issue that has drawn an enormous amount of attention over the last three years is the possible existence of statistical asymmetries in the behaviour of CMB radiation. The assumption of an isotropic universe is very much embedded at the roots of cosmology, in view of the so-called *Einstein cosmological principle* that the universe should "look the same" to any observer. However, quite unexpectedly, some evidence of statistical anisotropy is indeed present in the first releases of WMAP data (2003, 2006); see, for instance, Hansen *et al.* [10] and the references therein. Statistical procedures considered so far have led to inconclusive results. This mixed evidence has sparked an enormous amount of further empirical research, as such asymmetries may entail profound consequences in fundamental physics (as mentioned in the Introduction). It is therefore of great importance to devise new statistical tests that can exploit as efficiently as possible the available data. In view of their double localization in real and harmonic space, needlets emerge as natural candidates to build such procedures. We stress, in fact, that the knowledge of the scales where the asymmetries might lie would provide essential information toward their understanding.

We mention here that many other new wavelet-related systems have recently been introduced to deal with anisotropic features, such as curvelets and other forms of directional wavelets (see Jin *et al.* [13], Starck *et al.* [26], Vielva *et al.* [30] and Wiaux *et al.* [32]). We view these important approaches as complementary to ours and we leave for future research the investigation of interactions between these lines of research.

We suggest the following procedures: For fixed $q$ and $A$, consider

$$S_j := \sup_{a=1,\ldots,A_r} \left| \frac{\Gamma_{a;rj}}{\sqrt{\widehat{\Sigma}_j}} \right|;$$

under the null of Gaussianity and isotropy, we have from Theorem 11, part (c):

$$\lim_{j \to \infty} \Pr\{S_j \leq x\} = \{\Phi(x) - \Phi(-x)\}^{A_r},$$

where $\Phi$ denotes the standard cumulative distribution function of a Gaussian random variable. Under the alternative, we would expect to obtain unusually large values over some regions of the sky where data are generated according to a different model, and thus we expect the procedure to have good power properties against suitable alternatives.

An alternative approach can be envisaged as follows: In many applications, a possible direction for the asymmetries is easily conjectured. For instance, with CMB data, natural frames of reference are provided by either the ecliptic plane (i.e., the subspace approximately spanned by the planetary orbits around the sun) or the galactic plane, that is, the plane approximately defined by the location of the Milky Way. Asymmetries in these directions would lead presumably to explanations that do not directly involve the CMB itself, but rather other astrophysical entities of a non-cosmological nature. These forms of north–south asymmetries can be readily tested by means of the procedure that we describe below. Let

$$\Gamma_{ja;q} := \frac{1}{\sqrt{N_j}} \sum_{k \in \Xi_{\pi B^{-j}} \cap S(x_a)} H_q(\widehat{\beta}_{jk}), \qquad a = 1, 2,$$

where $x_1, x_2 \in \mathbb{S}^2$ denote the north and south poles in the suitable frame of reference (e.g., the so-called galactic or ecliptic poles). From Theorem 11, part (c), we have that

$$T_j := \frac{1}{2\widehat{\Sigma}_j} \{\Gamma_{j1;2} - \Gamma_{j2;2}\}^2 \to_d \chi_1^2 \qquad \text{under } H_0,$$

thus immediately providing threshold values with an asymptotic justification.

More sophisticated approaches are indeed possible. One idea is to exploit the asymptotically i.i.d. behaviour of needlet coefficients to transform these statistics into an approximate sample of spherical directional data, and then use the rich machinery of testing for uniformity that has been developed in these circumstances (see Mardia and Jupp [14], Pycke [24] and the references therein). A possible idea is to proceed with a hard thresholding of needlet coefficients and then focus on the directions that correspond to selected coefficients. More precisely, we can define the new data set

$$\{x_{rj;\tau} \in \mathbb{S}^2 : |\widehat{\beta}_{jk}| > \tau_j\},$$

for a suitable choice of the thresholding sequence $\{\tau_j\}$. Under the null, in view of the limiting properties of $\{\widehat{\beta}_{jk}\}$, the sequence $\{x_{rj;\tau}\}$ is approximately uniformly distributed on the sphere, a conclusion that can easily tested by a variety of well-established procedures. A rigorous investigation of this proposal, however, is beyond the scope of the present paper and will be deferred to future research. Likewise, applications of the ideas in this section to CMB data from WMAP are currently in progress and will be reported elsewhere.

## 6. Appendix: Proof of Proposition 13

For notational simplicity and without loss of generality, throughout this Appendix we take $B = 2$. The proof requires several steps. The first part, which is relatively simple, refers to circumstances where the distance between Voronoi cells is larger than their radius.

### 6.1. First case $d(S(x_a), S(x_b)) \geq \delta$

Under this hypothesis, and due to the next lemma, we have

$$\frac{1}{\sqrt{N_{a;\delta\varepsilon}N_{b;\delta\varepsilon}}} \sum_{v \in \Xi_\varepsilon \cap S(x_a)} \sum_{u \in \Xi_\varepsilon \cap S(x_b)} \frac{1}{(1 + 2^j d(u,v))^M}$$

$$\leq \sqrt{N_{a;\delta\varepsilon}N_{b;\delta\varepsilon}} \frac{1}{(2^j \delta)^M} \leq 2\pi^2 \left(\frac{\delta}{\varepsilon}\right)^2 \left(\frac{\varepsilon}{\delta}\right)^M = 2\pi^2 \left(\frac{\varepsilon}{\delta}\right)^{M-2}.$$

**Lemma 14.** *Let $0 < \varepsilon \leq \frac{\delta}{4}$ and let $\Xi_\delta$ be a maximal $\delta$-net and $\Xi_\varepsilon$ bea maximal $\varepsilon$-net. Let $x_a \in \Xi_\delta$ and let $S(x_a)$ be the corresponding Voronoi cell. Then*

$$\left(\frac{\delta}{\varepsilon}\right)^2 \frac{1}{4\pi^2} \leq \mathrm{Card}(S(x_a) \cap \Xi_\varepsilon) \leq 2\pi^2 \left(\frac{\delta}{\varepsilon}\right)^2. \tag{23}$$

**Proof.** As $B(x_a, \frac{\delta}{2}) \subset S(x_a) \subset B(x_a, \delta)$, we have

$$\bigcup_{u \in S(x_a) \cap \Xi_\varepsilon} B\left(u, \frac{\varepsilon}{2}\right) \subset B\left(x_a, \delta + \frac{\varepsilon}{2}\right)$$

and

$$B\left(x_a, \frac{\delta}{2} - \varepsilon\right) \subset \bigcup_{u \in S(x_a) \cap \Xi_\varepsilon} B(u, \varepsilon).$$

In view of (8), it follows easily that

$$\mathrm{Card}(\{u \in S(x_a) \cap \Xi_\varepsilon\})\frac{\varepsilon^2}{\pi} \leq \pi\left(\delta + \frac{\varepsilon}{2}\right)^2.$$

Hence

$$\mathrm{Card}(\{u \in S(x_a) \cap \Xi_\varepsilon)\}) \leq 2\pi^2\left(\frac{\delta}{\varepsilon}\right)^2$$

and moreover

$$\frac{4}{\pi}\left(\frac{\delta}{2} - \varepsilon\right)^2 \leq \mathrm{Card}(\{u \in S(x_a) \cap \Xi_\varepsilon\})\pi\varepsilon^2,$$

hence

$$\left(\frac{\delta}{\varepsilon}\right)^2 \frac{1}{4\pi^2} \leq \frac{4}{\pi^2}\left(\frac{\delta}{2\varepsilon} - 1\right)^2 \leq \mathrm{Card}(\{u \in S(x_a) \cap \Xi_\varepsilon\}). \qquad \square$$

## 6.2. Second case $d(S(x_a), S(x_b)) < \delta$

Here we face the situation where we have neighbouring Voronoi cells; the corresponding covariances are clearly harder to bound and we shall first introduce several lemmas.

**Lemma 15.** *Let $\varepsilon \sim 2^{-j}$ and $M \geq 3$. Let $x_a \neq x_b$ and $x_a, x_b \in \Xi_\delta$. Let $u \in \Xi_\varepsilon \cap S(x_a)$ be fixed. Then*

$$\sum_{v \in \Xi_\varepsilon \cap S(x_b)} \frac{1}{(1 + 2^j d(u, v))^M}$$

$$\leq C_M \frac{1}{(1 + 2^j d(u, S(x_b)))^{M-2}}.$$

$C_M$ *can be chosen equal to $2\pi 3^{2M-1}$.*

**Proof.** By straightforward manipulations, we obtain the bound

$$A(u) = \sum_{v \in \Xi_\varepsilon \cap S(x_b)} \frac{1}{(1 + 2^j d(u, v))^M}$$

$$= \sum_{v \in \Xi_\varepsilon \cap S(x_b)} \frac{1}{|B(v, \varepsilon/2)|} \int_{B(v, \varepsilon/2)} \frac{1}{(1 + 2^j d(u, v))^M} \, dx$$

$$\leq \frac{\pi}{\varepsilon^2} \sum_{v \in \Xi_\varepsilon \cap S(x_b)} \int_{B(v, \varepsilon/2)} \frac{2^M}{(1 + 2^j d(u, x))^M} \, dx$$

by (8) and as by the triangle inequality $d(u,x) \leq 2d(u,v)$ for $x \in B(v, \varepsilon/2)$. Clearly

$$\sum_{v \in \Xi_\varepsilon \cap S(x_b)} \int_{B(v,\varepsilon/2)} \frac{2^M}{(1+2^j d(u,x))^M} \, dx$$

$$\leq \int_{B^c(x_i,\varepsilon/2)} \frac{2^M}{(1+2^j d(u,x))^M} \, dx$$

$$\leq 2\pi \int_{\varepsilon/2}^\pi \frac{2^M \sin\theta}{(1+2^j \theta)^M} \, d\theta \leq 2\pi \int_{\varepsilon/2}^\pi \frac{2^M \theta}{(2^j \theta)^M} \, d\theta \leq \frac{\pi 2^{2M-1}}{M-2} 2^{-2j} \leq \pi 2^{2M-1} 2^{-2j},$$

which implies

$$A(u) \leq \pi 2^{2M-1}.$$

The previous inequality for $d(u, S(x_a)) \leq \frac{\varepsilon}{2}$ yields immediately

$$A(u) \leq \pi 3^{2M-1} \frac{1}{(1+2^j d(u, S(x_b)))^{M-2}}.$$

On the other hand, if $d(u, S(x_a)) > \frac{\varepsilon}{2}$,

$$\sum_{v \in \Xi_\varepsilon \cap S(x_b)} \int_{B(v,\varepsilon/2)} \frac{2^M}{(1+2^j d(u,x))^M} \, dx$$

$$\leq \int_{B^c(u,d(u,S(x_a)))} \frac{2^M}{(1+2^j d(u,x))^M} \, dx$$

$$\leq \int_{d(u,S(x_a))}^\pi \frac{2\pi 2^M \sin\theta}{(1+2^j \theta)^M} \, d\theta \leq 2\pi 2^{(1-j)M} \int_{d(u,S(x_a))}^\pi \theta^{1-M} \, d\theta$$

$$\leq 2\pi 2^{-j} \frac{2}{M-2} (2^j d(u, S(x_a)))^{-(M-2)} \leq 2\pi 2^{-j} (2^j d(u, S(x_a)))^{-(M-2)},$$

whence we get

$$A(u) \leq 2\pi (2^j d(u, S(x_b)))^{-(M-2)} \leq 2\pi 3^{M-2} \frac{1}{(1+2^j d(u, S(x_b)))^{M-2}}.$$

We are thus able to conclude that

$$\sum_{v \in \Xi\varepsilon \cap S(x_b)} \frac{1}{(1+2^j d(u,v))^M} \leq 2\pi 3^{2M-1} \frac{1}{(1+2^j d(u, S(x_b)))^{M-2}}. \qquad \square$$

Using the previous Lemma 15, we have:

$$W \stackrel{\text{def}}{=} \frac{1}{\sqrt{N_{a;\delta\varepsilon}N_{b;\delta\varepsilon}}} \sum_{u \in \Xi_\varepsilon \cap S(x_a)} \sum_{v \in \Xi_\varepsilon \cap S(x_b)} \frac{1}{(1+2^j d(u,v))^M}$$

$$\leq \frac{1}{\sqrt{N_{a;\delta\varepsilon}N_{b;\delta\varepsilon}}} 2\pi 3^{2M-1} \sum_{u \in \Xi_\varepsilon \cap S(x_a)} \frac{1}{(1+2^j d(u,S(x_b)))^{M-2}}.$$

Now as $d(S(x_a), S(x_b)) \leq \delta$, let us observe that if $u \in \Xi_\varepsilon \cap S(x_a)$

$$d(u, S(x_b)) \leq d(S(x_a), S(x_b)) + 2\delta \leq 3\delta.$$

Some computations yield, using Abel's formula of summation by parts:

$$W \leq \frac{2\pi 3^{2M-1}}{\sqrt{N_{a;\delta\varepsilon}N_{b;\delta\varepsilon}}} \sum_{0 \leq l \leq 3\delta/\varepsilon} \frac{\text{Card}\{u \in \Xi_\varepsilon \cap S(x_a), l\varepsilon \leq d(u,S(x_b)) < (l+1)\varepsilon\}}{(1+l)^{M-2}}$$

$$= \frac{2\pi 3^{2M-1}}{\sqrt{N_{a;\delta\varepsilon}N_{b;\delta\varepsilon}}} \sum_{1 \leq l \leq 3\delta/\varepsilon+1} \frac{\text{Card}\{u \in \Xi_\varepsilon \cap S(x_a), (l-1)\varepsilon \leq d(u,S(x_b)) < l\varepsilon\}}{l^{M-2}}$$

$$\leq \frac{2\pi 3^{2M-1}}{\sqrt{N_{a;\delta\varepsilon}N_{b;\delta\varepsilon}}} \Bigg\{ \left(\frac{\varepsilon}{3\delta}\right)^{M-2} \text{Card}\{u \in \Xi_\varepsilon \cap S(x_a), d(u,S(x_b)) < 3\delta + \varepsilon\}$$

$$+ \sum_{1 \leq l \leq 3\delta/\varepsilon} \left(\frac{1}{l^{M-2}} - \frac{1}{(1+l)^{M-2}}\right) \text{Card}\{u \in \Xi_\varepsilon \cap S(x_a), d(u,S(x_b)) < l\varepsilon\}\Bigg\}$$

$$\leq 4\pi^2 \left(\frac{\varepsilon}{\delta}\right)^2 2\pi 3^{2M-1} \Bigg\{ \left(\frac{\varepsilon}{3\delta}\right)^{M-2} 2\pi^2 \left(\frac{\delta}{\varepsilon}\right)^2$$

$$+ \sum_{1 \leq l \leq 3\delta/\varepsilon} \frac{M-2}{l^{M-1}} \text{Card}\{u \in \Xi_\varepsilon \cap S(x_a), d(u,S(x_b)) < l\varepsilon\}\Bigg\}$$

$$\leq 16\pi^5 3^{M+1} \left(\frac{\varepsilon}{\delta}\right)^{M-2}$$

$$+ 8\pi^3 3^{2M-1} \left(\frac{\varepsilon}{\delta}\right)^2 \sum_{1 \leq l \leq 3\delta/\varepsilon} \frac{M-2}{l^{M-1}} \text{Card}\{u \in \Xi_\varepsilon \cap S(x_a), d(u,S(x_b)) < l\varepsilon\}.$$

Let us introduce a further auxiliary result.

**Lemma 16.** *For $d(S(x_a), S(x_b)) \leq \delta$, $\varepsilon \leq \delta/4$,*

$$\text{Card}\{u \in \Xi_\varepsilon \cap S(x_a), d(u,S(x_b)) < l\varepsilon\} \leq C_1(l+1)^2 + C_2 \frac{\delta}{\varepsilon}(l+1)$$

*(with, e.g., $C_1 = 6\pi^4$ and $C_2 = 54\pi^4$).*

**Proof.** Using Lemma 8, if $\partial A$ denotes the boundary of the set $A$, we have

$$\partial(S(x_b)) = \bigcup_{m=1}^{M} \Gamma_m, \qquad M \leq 6\pi^2,$$

where $\Gamma_m = [a_m, a_{m+1}]$ is a portion of a geodesic circle $\Gamma_m \subset C_{c_m} := \partial(B(c_m, \frac{\pi}{2}))$ of length less than $\text{diam}(S(x_b)) \leq 2\delta$, $c_m$ denoting its center. Let $u \in \Xi_\varepsilon \cap S(x_a)$, $d(u, S(x_b)) < l\varepsilon$. Certainly

$$d(u, S(x_b)) = d(u, w), \qquad w \in \partial(S(x_b)).$$

We split the proof into two parts according to whether $w$ is a corner point or not.

- If $w$ is a corner point $a_m$, then

$$u \in \Xi_\varepsilon, \qquad d(u, a_m) \leq l\varepsilon \leq 3\delta$$

and

$$\text{Card}\{u \in \Xi_\varepsilon, d(u, a_m) \leq l\varepsilon\} \leq \tfrac{1}{4}\pi^2(2l+1).$$

In effect, for such $u$ we have

$$\bigcup_u B\left(u, \frac{\varepsilon}{2}\right) \subset B\left(a_m, \left(l + \frac{1}{2}\right)\varepsilon\right).$$

So if $k_m = \text{Card}\{u \in \Xi_\varepsilon, d(u, a_m) < l\varepsilon\}$, as $l \geq 1$, using (8)

$$k_m \frac{4}{\pi} \left(\frac{\varepsilon}{2}\right)^2 \leq \pi \left(l + \frac{1}{2}\right)^2 \varepsilon^2,$$

hence

$$k_m \leq \tfrac{1}{4}\pi^2(2l+1)^2 \leq \pi^2(l+1)^2.$$

- On the other hand, let us focus on $w \in ]a_m, a_{m+1}[$,

$$d(u, S(x_b)) = d(u, w) = d(u, C_{c_m}) \leq 3\delta.$$

So we have that

$$u \in B\left(c_m, \frac{\pi}{2}\right) \setminus B\left(c_m, \frac{\pi}{2} - l\varepsilon\right)$$

and also to the portion of the half sphere $B(c_m, \pi/2)$, which is between the two geodesics joining $c_m$ to $a_m$ and $c_m$ to $a_{m+1}$. Let us call this set $G(c_m; a_m, a_{m+1})$.

Therefore $\{u \in \Xi_\varepsilon \cap S(x_a), d(u, S(x_b)) < l\varepsilon\}$ is contained in the union of the $M$ sets $\{u \in \Xi_\varepsilon \cap B(a_m, l\varepsilon)\}$ and in the union of the $M$ sets $\{u \in \Xi_\varepsilon \cap G(c_m; a_m, a_{m+1}) \cap B(c_m, \frac{\pi}{2}) \setminus B(c_m, \frac{\pi}{2} - l\varepsilon)\}$. Let us evaluate the cardinality of each of these sets.

Let us consider the set $\{u \in \Xi_\varepsilon \cap G(c_m; a_m, a_{m+1}) \cap B(c_m, \frac{\pi}{2}) \setminus B(c_m, \frac{\pi}{2} - l\varepsilon)\}$. It is a simple observation that for a point $u$ of this set,

$$B\left(u, \frac{\varepsilon}{2}\right) \subset B\left(c_m, \frac{\pi}{2} + \frac{\varepsilon}{2}\right) \setminus B\left(c_m, \frac{\pi}{2} - \left(l + \frac{1}{2}\right)\varepsilon\right) \tag{24}$$

(recall that the $B(u, \frac{\varepsilon}{2})$ are disjoint), and also

$$B\left(u, \frac{\varepsilon}{2}\right) \subset G(c_m; a'_m, a'_{m+1}),$$

where $[a'_m, a'_{m+1}]$ is a portion of geodesic circle ($\Gamma_m \subset C_{c_m}$ containing $[a_m, a_{m+1}]$ of length $2\delta + 2\alpha$, with $\alpha \leq \delta$. Indeed, such a point $u$ is on a geodesic joining $c_m$ to $w \in ]a_m, a_{m+1}[$.

Let us recall that $d(c_m, w) = \frac{\pi}{2}$ and thus as $d(u, w) \leq 3\delta$ and if $\delta \leq \frac{\pi}{12}$, certainly, the ball $B(u, \varepsilon)$ does not contain the point $c_m$ (as $\varepsilon \leq \frac{\delta}{4}$).

Hence, let us consider the geodesic tangent $[c_m, v]$, where $v$ belongs to the boundary of the ball. Let us also consider the spherical triangle $(c_m, u, v)$. The $\alpha$ that we are trying to bound is the angle at the vertex $c_m$.

By a standard spherical trigonometric formula, we have:

$$\frac{\sin \alpha}{\sin d(u, v)} = \frac{\sin \hat{v}}{\sin d(u, c_m)},$$

where $\hat{v}$ denotes the angle at the vertex $v$ in the spherical triangle $(c_m, u, v)$. This formula translates here as:

$$\frac{\sin \alpha}{\sin \varepsilon} = \frac{\sin \pi/2}{\sin(\pi/2 - d(u, w))} = \frac{1}{\cos d(u, w)}.$$

Hence, as $\varepsilon \leq \frac{\delta}{4}$,

$$\sin \alpha \leq \frac{\sin \varepsilon}{\cos 3\delta} \leq \frac{\sin \delta/4}{\cos 3\delta} \leq \sin \delta.$$

The last inequality is easy to check for $\delta \leq \frac{\pi}{12}$.

We are now in the position to conclude the proof of Lemma 16.

Let us first take $\delta \leq \frac{\pi}{12}$. So, if $k_m = \text{Card}\{u \in \Xi_\varepsilon \cap G(c_m; a_m, a_{m+1}) \cap B(c_m, \frac{\pi}{2}) \setminus B(c_m, \frac{\pi}{2} - l\varepsilon)\}$. By (24)

$$k_m \frac{4}{\pi}\left(\frac{\varepsilon}{2}\right)^2 \leq \left| G(c_m; a_m, a_{m+1}) \cap B\left(c_m, \frac{\pi}{2} + \frac{\varepsilon}{2}\right) \setminus B\left(c_m, \frac{\pi}{2} - \left(l + \frac{1}{2}\right)\varepsilon\right) \right|.$$

So

$$k_m \frac{4}{\pi}\left(\frac{\varepsilon}{2}\right)^2 \leq \frac{4\delta}{2\pi}\pi\left[\left(\frac{\pi}{2}+\frac{\varepsilon}{2}\right)^2 - \left(\frac{\pi}{2}-\left(l+\frac{1}{2}\right)\varepsilon\right)^2\right]$$
$$= 2\delta\left[\pi(l+1)\varepsilon + \frac{1}{2}\varepsilon^2 - \left(l+\frac{1}{2}\right)\varepsilon^2\right] \leq 2\pi\delta(l+1)\varepsilon.$$

So

$$k_m \leq \frac{\delta}{\varepsilon}2\pi^2(l+1).$$

On the other hand, let us suppose instead that $\frac{\pi}{12} < \delta \leq \pi$. The computation of the cardinality $K = \mathrm{Card}\{u \in \Xi_\varepsilon \cap B(c_m, \frac{\pi}{2}) \setminus B(c_m, \frac{\pi}{2} - l\varepsilon)\}$ is much simpler: By (24)

$$K\frac{4}{\pi}\left(\frac{\varepsilon}{2}\right)^2 \leq \left|B\left(c_m, \frac{\pi}{2}+\frac{\varepsilon}{2}\right) \setminus B\left(c_m, \frac{\pi}{2}-\left(l+\frac{1}{2}\right)\varepsilon\right)\right|,$$

which leads to

$$K\frac{4}{\pi}\left(\frac{\varepsilon}{2}\right)^2 \leq \pi\left[\left(\frac{\pi}{2}+\frac{\varepsilon}{2}\right)^2 - \left(\frac{\pi}{2}-\left(l+\frac{1}{2}\right)\varepsilon\right)^2\right] \leq \pi^2(l+1)\varepsilon$$

and therefore

$$K \leq \frac{1}{\varepsilon}\pi^3(l+1) = \frac{\pi}{12}\frac{1}{\varepsilon}12\pi^2(l+1) \leq \frac{\delta}{\varepsilon}12\pi^2(l+1).$$

By the previous computations and recalling that the number of adjacent cells is $\leq 6\pi^2$, we have completed the proof of Lemma 16. □

Now we can conclude the proof of Proposition 13 by noting that

$$W \leq 16\pi^5 3^{M+1}\left(\frac{\varepsilon}{\delta}\right)^{M-2}$$
$$+ 8\pi^3 3^{M-2}\left(\frac{\varepsilon}{\delta}\right)^2 \sum_{1\leq l\leq 3\delta/\varepsilon} \frac{M-2}{l^{M-1}} \mathrm{Card}\{u \in \Xi_\varepsilon \cap S(x_i), d(u, S(x_b)) < l\varepsilon\}$$
$$\leq 16\pi^5\left(\frac{\varepsilon}{\delta}\right)^{M-2} + 8\pi^3 3^{M-2}\left(\frac{\varepsilon}{\delta}\right)^2 \sum_{1\leq l\leq 3\delta/\varepsilon} \frac{M-2}{l^{M-1}}C_1(l+1)^2$$
$$+ 8\pi^3 3^{M-2}\frac{\varepsilon}{\delta} \sum_{1\leq l\leq 3\delta/\varepsilon} \frac{M-2}{l^{M-1}}C_2(l+1),$$

and if $M \geq 3$,

$$\left(\frac{\varepsilon}{\delta}\right)^2 \sum_{1 \leq l \leq 3\delta/\varepsilon} \frac{1}{l^{M-1}}(l+1)^2 \sim \left(\frac{\varepsilon}{\delta}\right)^2 \sum_{1 \leq l \leq 3\delta/\varepsilon} \frac{1}{l^{M-3}} \leq C\frac{\varepsilon}{\delta}$$

as well as

$$\frac{\varepsilon}{\delta} \sum_{1 \leq l \leq 3\delta/\varepsilon} \frac{1}{l^{M-1}}(l+1) \sim \frac{\varepsilon}{\delta} \sum_{1 \leq l \leq 3\delta/\varepsilon} \frac{1}{l^{M-2}} \leq C\frac{\varepsilon}{\delta}\log\frac{\delta}{\varepsilon}.$$